\documentclass[a4paper,12pt,twoside]{article}
\usepackage{amsmath,amsthm,amssymb,yoshinori,cite}
\numberwithin{equation}{section}

\layout

\newcommand{\limup}{\lim_{q \uparrow 1}}

\newcommand{\nn}{^{(\n)}}
\newcommand{\tp}[1]{^{(#1)}}

\newcommand{\tB}{\widetilde{B}}

\newcommand{\bsym}[1]{{\boldsymbol{#1}}}

\newcommand{\f}[1]{F^{(#1)}_{q,0}}
\renewcommand{\fp}[1]{F^{(#1)}_{q,+}}
\renewcommand{\fm}[1]{F^{(#1)}_{q,-}}
\newcommand{\fpm}[1]{F^{(#1)}_{q,\pm}}
\renewcommand{\fe}[1]{F^{(#1)}_{q,\e}}

\renewcommand{\d}{\delta_q}

\newcommand{\ua}{\uparrow}

\boldtitle{Integral representations of $q$-analogues of\\ the Hurwitz
zeta function}  
\author{Masato Wakayama and Yoshinori Yamasaki}
\date{\today}

\begin{document}
%===========================================================================
%===========================================================================

\setlength{\baselineskip}{16pt}
\maketitle

\begin{abstract}
 Two integral representations of $q$-analogues of the Hurwitz zeta
 function are established. Each integral representation allows us to
 obtain an analytic continuation including also a full description of
 poles and special values at non-positive integers of the $q$-analogue
 of the Hurwitz zeta function, and to study the classical limit of this 
 $q$-analogue. All the discussion developed here is entirely different
 from the previous work in \cite{KawagoeWakayamaYamasaki}.

\smallbreak

 2000\ Mathematics Subject Classification:\, 11B68, 11M35

 Key words:\, Hurwitz's zeta function, Mellin transforms, $q$-analogue,
 Bernoulli polynomials, Poisson's summation formula
\end{abstract}

%===========================================================================
\section{Introduction}
\label{sec:intro}
%===========================================================================

 By the integral expression of the gamma function $\G(s)$, Hurwitz's
 zeta function $\z(s,z):=\sum^{\infty}_{n=0}(n+z)^{-s}$ is obtained by the 
 Mellin transform of the generating function $G(t,z)$ of Bernoulli
 polynomials $B_m(z)$ (see, e.g. \cite{WhittakerWatson1927}):  
\begin{align}
 \z(s,z)
&=\frac{1}{\G(s)}\int^{\infty}_{0}t^{s-1}G(t,z)\frac{dt}{t}
 \qquad (\Re(s)>1),\label{for:int Hurwitz}\\
 G(t,z):
&=\frac{te^{(1-z)t}}{e^t-1}=\sum^{\infty}_{m=0}(-1)^mB_m(z)\frac{t^m}{m!} \qquad(|t|<2\pi).\label{def:Gt}
\end{align}
 For any $0<a\le +\infty$ and $0<\e<\min\{a,2\pi\}$, $\z(s,z)$ is
 represented also as 
\begin{equation}
\label{for:cont Hurwitz}
 \z(s,z)
=\frac{\G(1-s)}{2\pi\I}\int_{C(\e,a)}\frac{(-t)^se^{(1-z)t}}{e^t-1}\frac{dt}{t}
+\frac{1}{\G(s)}\int^{\infty}_{a}t^{s-1}G(t,z)\frac{dt}{t},
\end{equation}
 where $C(\e,a)$ is a contour along the real axis from $a$ to $\e$,
 counterclockwise around the circle of radius $\e$ with center at the
 origin, and then along the real axis from $\e$ to $a$. This integral
 representation is straightforward from the idea due to Riemann in
 $1859$. We need this kind of segmentation at $t=a$ for handling
 delicate relations presented among several $q$-series in the
 discussion. Since the contour integral defines an entire function,
 \eqref{for:cont Hurwitz} provides a meromorphic continuation of
 $\z(s,z)$. Moreover, by the residue theorem, one shows that $\z(s,z)$ 
 has a simple pole at $s=1$ with residue $B_0(z)=1$ and
 $\z(1-m,z)=-B_m(z)/m$ $(m\in\bN)$. Furthermore, if we take $a=+\infty$, by 
 the residue theorem again, the contour integral \eqref{for:cont
 Hurwitz} yields the functional equation ($=$ the symmetric invariance
 for $s\leftrightarrow 1-s$) of the Riemann zeta function $\z(s)$
 \cite{Riemann1859}.   

 Let $0<q<1$ and $[z]_q:=(1-q^z)/(1-q)$ for $z\in\bC$. The following
 (Dirichlet-type) $q$-series has been studied in
 \cite{KawagoeWakayamaYamasaki}.
\[
 \z_q(s,t,z):=\sum^{\infty}_{n=0}\frac{q^{(n+z)t}}{{[n+z]_q}^s}
 \qquad (\Re(t)>0).
\]
 We put $\z\nn_q(s,z):=\z_q(s,s-\n,z)$ for $\n\in\bN$.
 The meromorphic continuation of $\z\nn_q(s,z)$ was obtained in two
 ways; one is by the binomial expansion, while the other is by the
 Euler-Maclaurin summation formula. Though the expression obtained by
 the binomial expansion has much advantage for describing the location
 of poles and special values at non-positive integers, it is difficult
 to determine whether $\z\nn_q(s,z)$ can give a proper $q$-analogue of
 $\z(s,z)$, i.e., if $\limup \z\nn_q(s,z)=\z(s,z)$ holds for any
 $s\in\bC$. Actually, the proof of the main assertion in
 \cite{KawagoeWakayamaYamasaki} which characterizes such proper
 $q$-analogues among the family of the functions $\z_q(s,\vp(s),z)$,
 $\vp(s)$ being a meromorphic function, could be achieved when
 employing the Euler-Maclaurin formula together with a careful piece of
 analysis. 

 The aim of the present paper is then, in contrast to the previous
 work, to study integral representations of $\z\nn_q(s,z)$ which are
 considered respectively as analogues of \eqref{for:int Hurwitz} and
 \eqref{for:cont Hurwitz}. Each integral representation allows us not
 only to see the aforementioned facts for $\z\nn_q(s,z)$ concerning the
 poles and special values but also to prove that $\z\nn_q(s,z)$ realizes 
 a proper $q$-analogue of $\z(s,z)$. During the course of a discussion
 for obtaining such integral representation, the Poisson summation
 formula plays an important role.  

 We remark that the works in \cite{Tsumura1999, Tsumura2001} have
 treated an integral representation of another $q$-analogue of Hurwitz's 
 zeta function (Mellin's transform). The present study is, however,
 different from those in \cite{Tsumura1999, Tsumura2001} in the sense
 that our $q$-analogue of the zeta function is given by exactly a
 Dirichlet-type $q$-series, whereas the $q$-analogue in
 \cite{Tsumura1999, Tsumura2001} has some extra term. Precisely, see
 Corollary\,2.4 \cite{KawagoeWakayamaYamasaki}.     

 The plan of the paper is as follows: In \Sec{integral}, we introduce
 functions $\fe{\n}(t,z)$ for $\e\in\{+,-,0\}$ and show that
 $\z\nn_q(s,z)$ is expressed as the Mellin transform of $\fp{\n}(t,z)$
 (\Prop{int}). In \Sec{main}, introducing a $q$-analogue $B\nn_m(z;q)$
 of the Bernoulli polynomial $B_m(z)$, we study relations between
 $\fm{\n}(t,z)$ and the generating function $G\nn_q(t,z)$ of
 $B\nn_m(z;q)$. Moreover, we show that $B\nn_m(z;q)\to B_m(z)$ as
 $q\ua 1$ for each $\n$ (\Thm{true Ber}). We then establish two integral 
 representations of $\z\nn_q(s,z)$ (\Thm{main}), and as a corollary, we
 show $\limup \z\nn_q(s,z)=\z(s,z)$ $(s\in\bC)$. The proof of these
 theorems and corollary are based on the result for $\f{\n}(t,z)$
 (\Prop{Poisson}) obtained by the Poisson summation formula. In the
 final section, we introduce some two-variable function $Z_q(s,t)$ by an 
 infinite product. If $s$ equals a positive integer $m$, $Z_q(m,t)$
 coincides with the inverse of Appell's $\cO$-function ($=$ the multiple 
 elliptic gamma function). We show that the logarithmic derivative of
 $Z_q(s,t)$ gives $\z_q(s,t):=\z_q(s,t,1)$ and obtain recurrence
 equations among $\z\nn_q(s):=\z\nn_q(s,1)$. 

 Throughout the paper, we assume $0<q<1$. The number $\n$ always
 represents a positive integer.

% 
%\[
% \xymatrix{
% & & *\txt{analytic\\continuation by} & \txt{papers} &
% *\txt{easy to see}&  \\ 
% & & *++[F--]{\txt{binomial theorem}}& 
% \txt{\cite{KanekoKurokawaWakayama2003}, \cite{KawagoeWakayamaYamasaki}}
% &\txt{location of poles\\special values} &  \\  
% & *++[][F]{\z\nn_q(s,z)} 
%\ar[ru]
%\ar[r]
%\ar[dr] 
% & *++[F=]{\txt{Mellin transform}} &
% *\txt{this article}&
% \txt{location of poles\\special values\\classical limit} \\ 
% & & *++[F--]{\txt{Euler-Maclaurin\\summation formula}} & 
% \txt{\cite{KanekoKurokawaWakayama2003}, \cite{KawagoeWakayamaYamasaki}} & \txt{classical limit}
% } \qquad
%\]
%

%===========================================================================
\section{An integral expression of $\z\nn_q(s,z)$}
\label{sec:integral}
%===========================================================================

 To obtain an integral expression of $\z\nn_q(s,z)$ as
 \eqref{for:cont Hurwitz} for $\z(s,z)$, we study functions
 $\fpm{\n}(t,z)$ defined as     
\begin{align}
  \fp{\n}(t,z):
&=t^{\n}\sum^{\infty}_{n=0}q^{-{\n}(n+z)}\exp\bigl(-tq^{-(n+z)}[n+z]_q\bigr)\nonumber\\
&=t^{\n}\sum^{\infty}_{n=0}q^{-{\n}(n+z)}\exp\bigl(t[-(n+z)]_q\bigr),\label{def:Fp}\\\fm{\n}(t,z):
&=t^{\n}\sum^{-1}_{n=-\infty}q^{-{\n}(n+z)}\exp\bigl(-tq^{-(n+z)}[n+z]_q\bigr)\nonumber\\
&=t^{\n}\sum^{\infty}_{n=1}q^{{\n}(n-z)}\exp\bigl(t[n-z]_q\bigr),\label{def:Fm}
\end{align}
 when
 $z\in D_q:=\bigl\{z\in\bC\,\bigl|\,|\Im(z)|<\frac{\pi}{2}\frac{1}{|\log{q}|}\bigr\}$.
 We first note the

%%%%%%%%%%%%%%%%%%%%%%%%%%%%%%%%%%%%%%%%%%%%%%%%%%%%%%%%%%%%%%%%%%%%%%%%%%%%
\begin{lem}
\label{lem:Fe conv}
 $(i)$\ Let $z\in D_q$. Put
 $R_q(z):=\bigl\{t\in\bC\,\bigl|\,|\arg(t)-(\Im(z))\log{q}|<\pi/2\bigr\}$.
 Here we assume $-\pi\le\arg(t)<\pi$. Then the function $\fp{\n}(t,z)$
 is holomorphic in $R_q(z)$.  

 $(ii)$\ The function $\fm{\n}(t,z)$ is entire. 

 $(iii)$\ The functions $\fpm{\n}(t,z)$ satisfy respectively 
\begin{equation}
\label{for:q-diff Fpm}
 \fpm{\n}(qt,z)
=e^{-t}\bigl(\fpm{\n}(t,z)\pm t^{\n}q^{{\n}(1-z)}e^{t[1-z]_q}\bigr).
\end{equation}
\end{lem}
\begin{proof}
 Since the condition $t\in R_q(z)$ implies 
 $\Re\bigl(t\exp(-\I \Im(z)\log{q})\bigr)>0$, it is easy to see that the 
 series $\fp{\n}(t,z)$ converges absolutely whenever $t\in R_q(z)$.
 Hence we have $(i)$. Since the function $[n-z]_q$ is bounded for
 $n\ge 0$, the series $\fm{\n}(t,z)$ converges absolutely for all
 $t\in\bC$, whence the assertion $(ii)$ follows. The functional
 equations \eqref{for:q-diff Fpm} are straightforward. 
\end{proof}

 Note that $R_q(z)\supset \bR_{>0}$, the positive real axis, for any
 $z\in D_q$. 

%%%%%%%%%%%%%%%%%%%%%%%%%%%%%%%%%%%%%%%%%%%%%%%%%%%%%%%%%%%%%%%%%%%%%%%%%%%%
\begin{lem}
\label{lem:Lebesgue}
 $(i)$\ For $a>0$ and $w>0$, $w^ae^{-w}\le (ae^{-1})^a$ holds.

 $(ii)$\ Put $z=x+\I y\in D_q$ $(x,y\in\bR)$ and
 $\b_y:=\cos(y\log{q})$. For $t>0$, we have 
\begin{equation}
\label{for:|F|}
 |\fp{\n}(t,z)|
\le \exp\Bigl(-t\frac{\b_y q^{-x}-1}{1-q}\Bigr)
\Bigl\{t^{\n}q^{-\n x}+\Bigl(\frac{\n e^{-1}}{\b_y}\Bigr)^{\n}
\frac{1}{1-\exp\bigl(-t\b_y q^{-x}\bigr)}\Bigr\}.
\end{equation}
 Further, suppose
 $z\in J_q:=\bigl\{z=x+\I y\in D_q\,\bigl|\, x>0,\ q^{-x}\cos(y\log{q})>1\bigr\}.$ Then the function $t^{\a}\fp{\n}(t,z)$
 is an integrable function on $[0,\infty)$ provided $\Re(\a)>0$.
\end{lem}
\begin{proof}
 The inequality $(i)$ is obvious. Using the relation 
 $[n+z]_q=1+q+\cdots+q^{n-1}+q^n[z]_q$, we have
\[
 \fp{\n}(t,z)
=\exp\bigl(-tq^{-z}[z]_q\bigr)
\Bigl\{t^{\n}q^{-\n z}+\sum^{\infty}_{n=1}\Bigl(tq^{-(n+z)}\Bigr)^{\n}\exp\bigl(-tq^{-(n+z)}\bigr)\prod^{n-1}_{j=1}\exp\bigl(-tq^{-{z-j}}\bigr)\Bigr\}. 
\]
 Note that $\Re(q^{-z})=\b_y q^{-x}$. Since 
 $q^{-{x-j}}\ge q^{-x}$ $(1\le j\le n-1)$, by $(i)$ with $a=\n$ and
 $w=t\b_y q^{-(n+x)}$, we get  
\begin{align*}
 |\fp{\n}(t,z)|
& \le \exp\Bigl(-t\frac{\b_y q^{-x}-1}{1-q}\Bigr)
\Bigl\{t^{\n}q^{-\n x}+\Bigl(\frac{{\n}e^{-1}}{\b_y}\Bigr)^{\n}
\sum^{\infty}_{n=1}\exp\bigl(-t{\b_y}q^{-x}(n-1)\bigr)\Bigr\}\\
& =\exp\Bigl(-t\frac{\b_y q^{-x}-1}{1-q}\Bigr)
\Bigl\{t^{\n}q^{-\n x}+\Bigl(\frac{\n e^{-1}}{\b_y}\Bigr)^{\n}
\frac{1}{1-\exp\bigl(-t\b_y q^{-x}\bigr)}\Bigr\}.
\end{align*} 
 This shows \eqref{for:|F|}. The rest of the assertion follows
 immediately from \eqref{for:|F|}. 
\end{proof}

 Estimate \eqref{for:|F|} verifies the expression of $\z\nn_q(s,z)$ by 
 the Mellin transform of $\fp{\n}(t,z)$.  

%%%%%%%%%%%%%%%%%%%%%%%%%%%%%%%%%%%%%%%%%%%%%%%%%%%%%%%%%%%%%%%%%%%%%%%%%%%%
\begin{lem}
\label{prop:int}
 Retain the notation in \Lem{Lebesgue}. Suppose $z\in J_q$. Then
\begin{equation}
\label{for:integral}
 \z\nn_q(s,z)
=\frac{1}{\G(s)}\int^{\infty}_{0}t^{s-\n}\fp{\n}(t,z)\frac{dt}{t} \qquad
(\Re(s)>\n+1). 
\end{equation}
\end{lem}
\begin{proof}
 Since $\Re(s)>\n+1$, $t^{s-\n-1}\fp{\n}(t,z)$ is integrable on $[0,\infty)$
 by \Lem{Lebesgue}. Hence we have  
\begin{align*}
\qquad\quad \int^{\infty}_{0}t^{s-\n}\fp{\n}(t,z)\frac{dt}{t}
&=\int^{\infty}_{0}t^s\sum^{\infty}_{n=0}q^{-{\n}(n+z)}\exp\bigl(-tq^{-(n+z)}[n+z]_q\bigr)\frac{dt}{t}\\
&=\sum^{\infty}_{n=0}\frac{q^{-{\n}(n+z)}}{\bigl(q^{-(n+z)}[n+z]_q\bigr)^s}\int^{\infty}_{0}w^s e^{-w}\frac{dw}{w}
   \qquad (w=tq^{-(n+z)}[n+z]_q)\\
&=\G(s)\z\nn_q(s,z).
\end{align*}
 The change of order of the integration and summation is legitimate by
 the Lebesgue convergence theorem. Hence the lemma follows.
\end{proof}

 Define the function $\f{\n}(t,z)$ by
\begin{align}
 \f{\n}(t,z):
&=\fp{\n}(t,z)+\fm{\n}(t,z)\nonumber\\
&=t^{\n}\sum_{n\in\bZ}q^{-{\n}(n+z)}\exp\bigl(-tq^{-(n+z)}[n+z]_q\bigr)
 \qquad (t\in R_q(z)).\label{def:F}
\end{align}
 It is clear that $\f{\n}(t,z)$ is periodic, that is,
 $\f{\n}(t,z)=\f{\n}(t,z+1)$. The following proposition is the key of
 our analysis and gives just the Fourier expansion of $\f{\n}(t,z)$.   

%%%%%%%%%%%%%%%%%%%%%%%%%%%%%%%%%%%%%%%%%%%%%%%%%%%%%%%%%%%%%%%%%%%%%%%%%%%%
\begin{prop}
\label{prop:Poisson}
 Let $z\in D_q$. For $t\in R_q(z)$, we have
\begin{equation}
\label{for:Poisson}
 \f{\n}(t,z)
=-\frac{(1-q)^{\n}}{\log{q}}e^{\frac{t}{1-q}}\sum_{m\in\bZ}\Bigl(\frac{1-q}{t}\Bigr)^{m\d}\G({\n}+m\d)e^{2\pi\I mz},
\end{equation}
 where $\d:=2\pi\I/\log{q}$.
\end{prop}
\begin{proof}
 Let
 $f\nn_q(\eta,z):=q^{-{\n}(\eta+z)}\exp\bigl(-tq^{-(\eta+z)}[\eta+z]_q\bigr)$.   Put $a=t/(1-q)$. Then the Fourier transform $\hat{f}\nn_q(\x,z)$ of
 $f\nn_q(\eta,z)$ is calculated as  
\begin{align*}
\qquad\qquad \hat{f}\nn_q(\x,z)
&=\int^{\infty}_{-\infty}f\nn_q(\eta,z)e^{-2\pi\I \eta\x}d{\eta}\\
&=-\frac{e^{a+2\pi\I {\x}z}}{\log{q}}\int^{\infty}_{-\infty}e^{-{\n}w}e^{-ae^{-w}}e^{-\d w\x}dw &(w=(\eta+z)\log{q})\\
&=-\frac{e^{a+2\pi\I {\x}z}}{\log{q}}\int^{\infty}_{0}{\a^{\n}}e^{-a\a}{\a}^{\d\x}\frac{d\a}{\a}  &(\a=e^{-w})\qquad\quad\ \, \\
&=-\frac{e^{a+2\pi\I {\x}z}}{\log{q}}\Bigl(\frac{1-q}{t}\Bigr)^{\n+\d\x}\G({\n}+\d\x).
\end{align*}
 By the Poisson summation formula, we have
\[
 \f{\n}(t,z)
=t^{\n}\sum_{n\in\bZ}f\nn_q(n,z)
=t^{\n}\sum_{m\in\bZ}\hat{f}\nn_q(m,z).
\]
 Hence the desired formula follows.  
\end{proof}

\smallbreak
\smallbreak
 
 Let $0<a<+\infty$ and $m\in\bZ$. To obtain an analytic continuation of
 $\z\nn_q(s,z)$ as a function of $s$ via \eqref{for:integral}, it is
 useful to define the function $\vp\nn_m(s;a,q)$ as   
\begin{equation}
\label{def:vp}
 \vp\nn_m(s;a,q):=\int^{a}_{0}t^{s-\n-1-m\d}e^{\frac{t}{1-q}}dt.
\end{equation}
 Noting $|\vp\nn_m(s;a,q)|\le \vp\nn_0(\Re(s);a,q)$, we see that
 $\vp\nn_m(s;a,q)$ is holomorphic in $\Re(s)>\n$. 

%%%%%%%%%%%%%%%%%%%%%%%%%%%%%%%%%%%%%%%%%%%%%%%%%%%%%%%%%%%%%%%%%%%%%%%%%%%%
\begin{prop}
\label{prop:vp}
 The function $\vp\nn_m(s;a,q)$ admits a meromorphic continuation to the  
 whole plane. It has simple poles at $s=n+m\d$
 $(n\in\bZ_{\le \n})$ with 
\begin{equation}
\label{for:residue vp}
 \Res{s=n+m\d}\vp\nn_m(s;a,q)=\frac{1}{(\n-n)!}\frac{1}{(1-q)^{-n+\n}}.
\end{equation}
\end{prop}
\begin{proof}
 Assume $\Re(s)>\n$. For $l=-n\in\bZ_{\ge -\n}$, integration by parts
 yields
\begin{multline*}
\qquad\qquad \vp\nn_m(s;a,q)
=\sum^{l+\n+1}_{j=1}\frac{1}{(s-\n-m\d)_j}\frac{(-1)^{j-1}}{(1-q)^{j-1}}a^{s-\n-m\d+j-1}e^{\frac{a}{1-q}}\\  
+\frac{1}{(s-\n-m\d)_{l+\n+1}}\frac{(-1)^{l+\n+1}}{(1-q)^{l+\n+1}}\vp\nn_m(s+l+\n+1;a,q).\qquad\qquad\qquad\qquad
\end{multline*}
 Here $(s)_k:=s(s+1)\cdots (s+k-1)$ for $k\in\bN$. This gives an
 analytic continuation of $\vp\nn_m(s;a,q)$ to the region 
 $\Re(s)>-l-1$. Moreover, since 
\[
 \vp\nn_m(s+l+\n+1;a,q)\Bigl|_{s=-l+m\d}
=\int^{a}_{0}e^{\frac{t}{1-q}}dt
=(1-q)\bigl(e^{\frac{a}{1-q}}-1\bigr),
\]
 we have \eqref{for:residue vp}. 
\end{proof}

%===========================================================================
\section{$q$-Analogue of Bernoulli polynomials and the main theorem} 
\label{sec:main}
%===========================================================================

 First we study a $q$-analogue $B\nn_m(z;q)$ of the Bernoulli polynomial 
 $B_m(z)$. Consider the $q$-difference equation      
\begin{align}
\label{for:q-diff G}
 \left\{
\begin{array}{l}
\ \displaystyle{G\nn_q(0,z)}
=\displaystyle{-\frac{(1-q)^{\n}}{\log{q}}(\n-1)!}, \\[10pt]
\displaystyle{G\nn_q(qt,z)}
=\displaystyle{e^{-t}\bigl(G\tp{\n}_q(t,z)+t^{\n}q^{{\n}(1-z)}e^{t[1-z]_q}\bigr)}.
\end{array}
\right.
\end{align}

%%%%%%%%%%%%%%%%%%%%%%%%%%%%%%%%%%%%%%%%%%%%%%%%%%%%%%%%%%%%%%%%%%%%%%%%%%%%
\begin{lem}
\label{lem:FG}
 A solution of \eqref{for:q-diff G} which is continuous at $t=0$ is 
 $(uniquly)$ given by 
\begin{equation}
\label{for:FG}
 G\nn_q(t,z)
=-\frac{(1-q)^{\n}}{\log{q}}(\n-1)!e^{\frac{t}{1-q}}-\fm{\n}(t,z).
\end{equation}
 In particular, $G\nn_q(t,z)$ is holomorphic at $t=0$.
\end{lem}
\begin{proof}
 Put $H\nn_q(t,z):=e^{-\frac{t}{1-q}}G\nn_q(t,z)$. Then $H\nn_q(t,z)$
 satisfies
\begin{align}
\label{for:q-diff H}
\left\{
\begin{array}{l}
\ \displaystyle{H\nn_q(0,z)}
=\displaystyle{-\frac{(1-q)^{\n}}{\log{q}}(\n-1)!},\\[10pt]
\displaystyle{H\nn_q(qt,z)}
=\displaystyle{H\tp{\n}_q(t,z)+t^{\n}q^{{\n}(1-z)}\exp\Bigl(-\frac{q^{1-z}}{1-q}t\Bigr)}.
\end{array}
\right.
\end{align} 
 It follows that
\[
  H\nn_q(q^mt,z)-H\nn_q(t,z)
=t^{\n}\sum^{m}_{n=1}q^{{\n}(n-z)}\exp\Bigl(-\frac{q^{n-z}}{1-q}t\Bigr) \qquad (m\in\bN).
\] 
 Since the function $H\nn_q(t,z)$ is continuous at $t=0$, letting
 $m\to\infty$, we have
\begin{equation}
\label{for:H}
 H\nn_q(t,z)
=-\frac{(1-q)^{\n}}{\log{q}}(\n-1)!-t^{\n}\sum^{\infty}_{n=1}q^{{\n}(n-z)}\exp\Bigl(-\frac{q^{n-z}}{1-q}t\Bigr).
\end{equation}
 This proves \eqref{for:FG}. 
\end{proof}

%%%%%%%%%%%%%%%%%%%%%%%%%%%%%%%%%%%%%%%%%%%%%%%%%%%%%%%%%%%%%%%%%%%%%%%%%%%%
\begin{remark}
 We give another proof of \eqref{for:FG} using the functional
 equations \eqref{for:q-diff Fpm} and \eqref{for:q-diff G}. Let
 $\psi\nn_q(t,z):=G\nn_q(t,z)+\fm{\n}(t,z)$. Since $\fm{\n}(0,z)=0$, we 
 have $\psi\nn_q(0,z)=G\nn_q(0,z)$. Also, since
 $\psi\nn_q(qt,z)=e^{-t}\psi\nn_q(t,z)$, we see that 
 $\psi\nn_q(q^mt,z)=e^{-[m]_qt}\psi\nn_q(t,z)$ for $m\in\bN$. Hence by
 the continuity of $\psi\nn_q(t,z)$ at $t=0$, we have
 $\psi\nn_q(t,z)=\lim_{m\to \infty}e^{[m]_qt}\psi\nn_q(q^mt,z)=e^{\frac{t}{1-q}}\psi\nn_q(0,z)$. This shows \eqref{for:FG}.
\end{remark}

 By the Taylor expansions of the exponential functions in
 \eqref{for:FG} and \eqref{for:H}, we may write    
\begin{equation}
\label{def:tB}
 G\nn_q(t,z)=\sum^{\infty}_{m=0}(-1)^m\tB\nn_m(z;q)\frac{t^m}{m!} \qquad
 (t\in\bC),
\end{equation}
 i.e., the series expansion has infinite radius of convergence at $t=0$
 when $0<q<1$. Then we see from \eqref{for:q-diff G} that the function
 $\tB\nn_m(z;q)$ satisfies the recursion formula:       
\begin{align}
\label{for:recursion-tBn}
\left\{
\begin{array}{l}
 \displaystyle{\tB\nn_0(z;q)}
=\displaystyle{-\frac{(1-q)^{\n}}{\log{q}}(\n-1)!}, \\[10pt]
 \displaystyle{\sum^n_{m=0}(-1)^m\binom{n}{m}q^m\tB\nn_m(z;q)}
=\displaystyle{(-1)^n\tB\tp{\n}_n(z;q)+{\n}!\binom{n}{\n}q^{{\n}(1-z)}{[1-z]_q}^{n-\n}}.
\end{array}
\right.
\end{align}
 Note that, in particular,
\begin{equation}
\label{for:Btil}
 \tB\nn_n(z;q)=(-1)^{n+1}\frac{(1-q)^{\n-n}}{\log{q}}(\n-1)! \qquad
 (0\le n\le \n-1).
\end{equation}
 Define also the functions $\{B\nn_m(z;q)\}_{m\ge 0}$ by 
\begin{equation}
\label{def:B}
 B\nn_m(z;q)
:=(-1)^{\n-1}\frac{m!}{(m+\n-1)!}\tB\nn_{m+\n-1}(z;q).
\end{equation}
 It is clear that $B\tp{1}_m(z;q)=\tB\tp{1}_m(z;q)$. The first  
 few are given by  
\begin{align*}
 B\nn_0(z;q)&=\frac{q-1}{\log{q}}, \quad
  B\nn_1(z;q)=-\frac{q^{-{\n}z}}{1-q^{-\n}}+\frac{1}{\n\log{q}},\\
 B\nn_2(z;q)&=\frac{2}{q-1}\Bigl(-\frac{q^{-{\n}z}}{1-q^{-\n}}+\frac{q^{z(-\n-1)}}{1-q^{-\n-1}}+\frac{1}{\n(\n+1)\log{q}}\Bigr), \quad \ldots.
\end{align*}
 From \Lem{FG} and \eqref{def:B}, we have a closed expression of
 $B\nn_m(z;q)$. 

%%%%%%%%%%%%%%%%%%%%%%%%%%%%%%%%%%%%%%%%%%%%%%%%%%%%%%%%%%%%%%%%%%%%%%%%%%%%
\begin{prop}
 For $m\in\bZ_{\ge 0}$,
\begin{equation}
\label{for:closed B_q}
  B\nn_m(z;q)
=(q-1)^{-m+1}\Biggl\{\sum^{m}_{l=1}(-1)^l\binom{m}{l}\frac{lq^{z(-l-\n+1)}}{1-q^{-l-\n+1}}+\binom{m+\n-1}{\n-1}^{-1}\frac{1}{\log{q}}\Biggr\}\in\bC[q^{-z}].
\end{equation}
\end{prop}\qed

 The following result shows the classical limit of $B\nn_m(z;q)$
 reproduces the Bernoulli polynomial. 

%%%%%%%%%%%%%%%%%%%%%%%%%%%%%%%%%%%%%%%%%%%%%%%%%%%%%%%%%%%%%%%%%%%%%%%%%%%%
\begin{thm}
\label{thm:true Ber}
 For $0<t<2\pi$, we have
\begin{equation}
\label{for:lim G}
 \limup G\nn_q(t,z)=t^{\n-1}G(t,z).
\end{equation}
 In particular, $\limup B\nn_m(z;q)=B_m(z)$ for all $m\in\bZ_{\ge 0}$. 
\end{thm}
\begin{proof}
 For any $t\in R_q(z)$, it is obvious that
\begin{equation}
\label{for:lim Fp}
 \limup \fp{\n}(t,z)
=t^{\n}\sum^{\infty}_{n=0}e^{-t(n+z)}
=t^{\n-1}\frac{te^{(1-z)t}}{e^t-1}
=t^{\n-1}G(t,z).
\end{equation}
 Also, from \eqref{for:Poisson} and \eqref{for:FG}, we have
\begin{equation}
\label{for:onR}
 G\nn_q(t,z)
=\frac{(1-q)^{\n}}{\log{q}}e^{\frac{t}{1-q}}\sum_{m\in\bZ\bslo}\Bigl(\frac{1-q}{t}\Bigr)^{m\d}\G({\n}+m\d)e^{2\pi\I mz}+\fp{\n}(t,z).
\end{equation}
 Then the series on the right hand side of \eqref{for:onR}
 disappears when $q\ua 1$ by the Stirling formula (see
 \cite{AndrewsAskeyRoy1999});   
\begin{equation}
\label{for:stirling}
  |\G(\n+m\d)|
\sim  \frac{(2\pi)^{\n}|m|^{\n-\frac{1}{2}}}{|\log{q}|^{\n-\frac{1}{2}}}e^{-\frac{{\pi}^2|m|}{|\log{q}|}}
\qquad (q\ua 1).
\end{equation}
 In fact, since $0<t<2\pi<\frac{3}{4}{\pi}^2$ and
 $1/\log{q}=-1/(1-q)+1/2+O(1-q)$, it follows that
\begin{align*}
 \frac{(1-q)^{\n}}{\log{q}}\bigl|e^{\frac{t}{1-q}}\G(\n+m\d)\bigr|
& \sim
 (-1)^{\n-\frac{1}{2}}(2\pi)^{\n}|m|^{\n-\frac{1}{2}}\Bigl(\frac{1-q}{\log{q}}\Bigr)^{\n+\frac{1}{2}}\frac{e^{-\frac{1}{4}\frac{{\pi}^2|m|}{|\log{q}|}}}{(1-q)^{\frac{1}{2}}}e^{\frac{t}{1-q}}e^{-\frac{3}{4}\frac{{\pi}^2|m|}{|\log{q}|}}\\
& =O\Bigl(\frac{e^{-\frac{1}{4}\frac{{\pi}^2|m|}{|\log{q}|}}}{(1-q)^{\frac{1}{2}}}\Bigr)\cdot \exp\Bigl(-\frac{\frac{3}{4}{\pi}^2|m|-t}{1-q}+O(1)\Bigr) \to 0 \qquad (m\ne 0)
\end{align*}
 when $q\ua 1$. Hence, letting $q\ua 1$ in \eqref{for:onR}, we obtain
 \eqref{for:lim G} from \eqref{for:lim Fp}. This completes the proof.
\end{proof}

%%%%%%%%%%%%%%%%%%%%%%%%%%%%%%%%%%%%%%%%%%%%%%%%%%%%%%%%%%%%%%%%%%%%%%%%%%%%
\begin{remark}
 The function $B\tp{1}_m(z;q)$ appeared essentially in
 \cite{Tsumura1999, Tsumura2001} (see also \cite{Koblitz1982} and
 \cite{Satoh1989}). Actually, one can show that
 $B\tp{1}_m(z;q):=(-1)^m\widetilde{\b}_m(1-z,q)$, where
 $\widetilde{\b}_m(z,q)$ is the one in \cite{Tsumura2001}. Also, in 
 \cite{KanekoKurokawaWakayama2003}, by the explicit expression
 \eqref{for:closed B_q} for $z=\n=1$, it was shown in a different way
 that $\limup B\tp{1}_m(1;q)=B_m$ $(m\in\bZ_{\ge 0})$, where
 $B_m=B_m(1)$ is the Bernoulli number.
\end{remark}

 The following theorem is the main result.

%%%%%%%%%%%%%%%%%%%%%%%%%%%%%%%%%%%%%%%%%%%%%%%%%%%%%%%%%%%%%%%%%%%%%%%%%%%% 
\begin{thm}
\label{thm:main}
 Let $z\in J_q=\bigl\{z=x+\I y\in \bC\,\bigl|\, x>0,\ |y|<\frac{\pi}{2}\frac{1}{|\log{q}|},\ q^{-x}\cos(y\log{q})>1\bigr\}$. Assume $\Re(s)>\n+1$. Then
\begin{enumerate}
 \item For any $0<a<+\infty$ and $0<\e\le a$, we have
 \begin{multline}
\label{for:another}
 \z\nn_q(s,z)
=\frac{(-1)^{\n}\G(1-s)}{2\pi\I}\int_{C(\e,a)}(-t)^{s-\n}G\nn_q(t,z)\frac{dt}{t}+\frac{1}{\G(s)}\int^{\infty}_{a}t^{s-\n}\fp{\n}(t,z)\frac{dt}{t}\\
-\frac{1}{\G(s)}\cdot \frac{(1-q)^{\n}}{\log{q}}\sum_{m\in\bZ\bslo}(1-q)^{m\d}\G(\n+m\d)e^{2\pi\I mz}\vp\nn_m(s;a,q),
\end{multline}
 where $C(\e,a)$ is the same contour as the one in \eqref{for:cont Hurwitz}.
 \item For any $N\in\bN$, we have  
\begin{multline}
\label{for:another2}
 \z\nn_q(s,z)
=\frac{1}{\G(s)}\Biggl\{\int^{1}_{0}t^{s-\n-1}\Bigl(G\nn_q(t,z)-\sum^{N+\n-1}_{k=0}(-1)^k\tB\nn_k(z;q)\frac{t^k}{k!}\Bigr)dt\\
+\sum^{N+\n-1}_{k=0}(-1)^k\frac{\tB\nn_k(z;q)}{k!}\frac{1}{s-\n+k}+\int^{\infty}_{1}t^{s-\n}\fp{\n}(t,z)\frac{dt}{t}\\
-\frac{(1-q)^{\n}}{\log{q}}\sum_{m\in\bZ\bslo}(1-q)^{m\d}\G(\n+m\d)e^{2\pi\I mz}\vp\nn_m(s;q)\Biggr\},
\end{multline}
 where $\vp\nn_m(s;q):=\vp\nn_m(s;1,q)$.
\end{enumerate}
\noindent
 Each integral representation shows that $\z\nn_q(s,z)$ is meromorphic
 in $\bC$ and has simple poles at the points in $n+\d\bZ$
 $(1\le n\le \n)$ and $\bZ_{\le 0}+\d\bZ\bslo$ with   
\begin{equation}
\label{for:residue zq}
\qquad \Res{s=n+m\d}\z\nn_q(s,z)
=-\binom{\n-1+m\d}{\n-n}\frac{(1-q)^{n+m\d}}{\log{q}}e^{2\pi\I mz}
\quad (n\in\bZ_{\le \n},\ m\in\bZ).
\end{equation}
 These exhaust all poles of $\z\nn_q(s,z)$. Further,
\begin{equation}
\label{for:sv-znq}
 \z\nn_q(1-m,z)=-\frac{B\nn_m(z;q)}{m} \qquad (m\in\bZ_{\ge 0}).
\end{equation}
\end{thm}
\begin{proof}
 For any $a>0$, from \eqref{for:integral}, we have  
\begin{equation}
\label{for:half}
 \G(s)\z\nn_q(s,z)
=\int^{a}_{0}t^{s-\n}\fp{\n}(t,z)\frac{dt}{t}+\int^{\infty}_{a}t^{s-\n}\fp{\n}(t,z)\frac{dt}{t}.
\end{equation}
 The second integral obviously defines an entire function by \Lem{Lebesgue}. From \eqref{for:onR} and
 \eqref{for:stirling} again, the first one turns to be 
\begin{equation}
\label{for:middle}
 \int^{a}_{0}t^{s-\n}G\nn_q(t,z)\frac{dt}{t}-\frac{(1-q)^{\n}}{\log{q}}\sum_{m\in\bZ\bslo}(1-q)^{m\d}\G(\n+m\d)e^{2\pi\I mz}\vp\nn_m(s;a,q).
\end{equation}
\begin{enumerate}
 \item Define the function $I\nn_q(s,z;a)$ by  
\[
 I\nn_q(s,z;a):
=\int_{C(\e,a)}(-t)^{s-\n}G\nn_q(t,z)\frac{dt}{t} \qquad (0<\e\le a). 
\]
       The integral converges absolutely and uniformly with respect to
       $s$ and does not depend on $\e$ by the Cauchy integral
       theorem. Hence $I\nn_q(s,z;a)$ is entire as a function of
       $s$. Moreover, we have   
\begin{equation}
\label{for:cont}
 \int^{a}_{0}t^{s-\n}G\nn_q(t,z)\frac{dt}{t}
=\frac{(-1)^{\n}\G(s)\G(1-s)}{2\pi\I}I\nn_q(s,z;a) \qquad (\Re(s)>\n+1).
\end{equation} 
       In fact, it is easy to see that 
\begin{align}
\label{for:iinntt}
 I\nn_q(s,z;a)
&=\Biggl(\int^{\e}_{a}+\int_{|t|=\e}+\int^{a}_{\e}\Biggr)(-t)^{s-\n}G\nn_q(t,z)\frac{dt}{t}\nonumber\\
&=(-1)^{\n}(e^{\pi s\I}-e^{-\pi s\I})\int^{a}_{\e}t^{s-\n}G\nn_q(t,z)\frac{dt}{t}\nonumber\\
&\ \ \ +\int_{|t|=\e}(-t)^{s-\n}G\nn_q(t,z)\frac{dt}{t}.
\end{align}
       Since $G\nn_q(t,z)$ is bounded for $t=|\e|$, the last integral in
       \eqref{for:iinntt} disappears when $\e\to 0$ provided
       $\Re(s)>\n+1$. Therefore, letting $\e\to 0$ in \eqref{for:iinntt}
       and using the relation $\G(s)\G(1-s)=\pi/\sin(\pi s)$, we obtain
       \eqref{for:cont}, whence \eqref{for:another}. Since
       $I\nn_q(s,z;a)$ is entire, \eqref{for:another} gives a
       meromorphic continuation of $\z\nn_q(s,z)$ to the whole plane
       $\bC$. By the residue theorem, it follows from \eqref{def:tB} that
\begin{equation}
\label{for:In}
 I\nn_q(n,z;a)=
\begin{cases}
 \displaystyle{\frac{2\pi\I}{(\n-n)!}\tB\nn_{\n-n}(z;q)}
 &\textrm{if}\quad n=1,2,\ldots,\n,\\
 \ \ \ \ 0 & \textrm{if}\quad n\in\bZ_{\ge \n+1}.
\end{cases}
\end{equation}
       Hence, by \Prop{vp}, we see that all poles of
       $\z\nn_q(s,z)$ are simple and given by $n+\d\bZ$ $(1\le n\le \n)$
       and $\bZ_{\le 0}+\d\bZ\bslo$. Further, if $s=n$ for
       $1\le n\le \n$, we have from \eqref{for:Btil} and \eqref{for:In}
       that  
\begin{align*}
 \Res{s=n}\z\nn_q(s,z)
&=(-1)^{\n-n}\tB\nn_{\n-n}(z;q)\frac{1}{(n-1)!(\n-n)!}
=-\binom{\n-1}{\n-n}\frac{(1-q)^n}{\log{q}}.
\intertext{If $s=n+m\d$ with $m\ne 0$, we have from
 \eqref{for:residue vp}}  
 \Res{s=n+m\d}\z\nn_q(s,z)
&=-\frac{(1-q)^{n+m\d}e^{2\pi\I mz}}{\log{q}}\frac{\G(\n+m\d)}{\G(n+m\d)(\n-n)!}\\
&=-\binom{\n-1+m\d}{\n-n}\frac{(1-q)^{n+m\d}}{\log{q}}e^{2\pi\I mz}.
\end{align*}
       Hence \eqref{for:residue zq} follows. From \eqref{for:cont}
       again, it follows that  
\begin{align*}
  \z\nn_q(1-m,z)
&=\frac{(-1)^{\n}\G(m)}{2\pi\I}I\nn_q(1-m,z;a)\\
&=-(-1)^{\n-1}\frac{(m-1)!}{(m+\n-1)!}\tB\nn_{m+\n-1}(z;q)
=-\frac{B\nn_m(z;q)}{m}.
\end{align*}

 \item  Since the first integral in \eqref{for:another2} converges
       absolutely for $\Re(s)>1-N$, it suffices to show
       \eqref{for:another2}. This is easy because the integral
       $\int^{1}_{0}t^{s-\n-1}G\nn_q(t,z)dt$ in \eqref{for:middle} can
       be written as 
\begin{multline*}
\label{for:midddle}
 \int^{1}_{0}t^{s-\n-1}\Bigl(G\nn_q(t,z)-\sum^{N+\n-1}_{k=0}(-1)^k\tB\nn_k(z;q)\frac{t^k}{k!}\Bigr)dt
+\int^{1}_{0}t^{s-\n-1}\sum^{N+\n-1}_{k=0}(-1)^k\tB\nn_k(z;q)\frac{t^k}{k!}dt\\
=\int^{1}_{0}t^{s-\n-1}\Bigl(G\nn_q(t,z)-\sum^{N+\n-1}_{k=0}(-1)^k\tB\nn_k(z;q)\frac{t^k}{k!}\Bigr)dt
+\sum^{N+\n-1}_{k=0}(-1)^k\frac{\tB\nn_k(z;q)}{k!}\frac{1}{s-\n+k}.
\end{multline*}
\end{enumerate}
 This completes the proof of the theorem.
\end{proof}

%%%%%%%%%%%%%%%%%%%%%%%%%%%%%%%%%%%%%%%%%%%%%%%%%%%%%%%%%%%%%%%%%%%%%%%%%%%%
\begin{remark}
 The binomial theorem yields the following series expression of
 $\z\nn_q(s,z)$ (see \cite{KawagoeWakayamaYamasaki}, also
 \cite{KanekoKurokawaWakayama2003}), which shows that $\z\nn_q(s,z)$ is
 meromorphic in $\bC$:
\[
  \z\nn_q(s,z)
=(1-q)^s\sum^{\infty}_{l=0}\binom{s+l-1}{l}\frac{q^{z(s-\n+l)}}{1-q^{s-\n+l}}.
\]
 Using the expression, we can get the information about poles and
 special values.
\end{remark}

 By \eqref{for:residue zq}, we note that  
\[
 \limup \Res{s=n}\z\nn_q(s,z)=
\begin{cases}
 1 & \textrm{if}\quad n=1,\\
 0 & \textrm{if}\quad 2\le n\le \n.
\end{cases}
\]
 This shows the poles of $\z\nn_q(s,z)$ at $s=2,3,\ldots,\n$ disappear
 when $q\ua 1$. We now see the classical limit of $\z\nn_q(s,z)$. In
 fact,   

%%%%%%%%%%%%%%%%%%%%%%%%%%%%%%%%%%%%%%%%%%%%%%%%%%%%%%%%%%%%%%%%%%%%%%%%%%%%
\begin{cor}
\label{cor:lim zq}
 We have
\begin{equation}
\label{for:pre-main}
 \limup \z\nn_q(s,z)=\z(s,z) \qquad (s\in\bC,\ s\ne1,2,\ldots,\n).
\end{equation}
\end{cor}
\begin{proof}
 By comparing formulas \eqref{for:cont Hurwitz} and \eqref{for:another}
 when $a=1$ (one may take any $a$ satisfying $0<a<2\pi$), in view of
 \Thm{true Ber}, it is sufficient to show  
\begin{equation}
\label{for:lim-main}
 \limup \frac{(1-q)^{\n}}{\log{q}}\sum_{m\in\bZ\bslo}(1-q)^{m\d}\G(\n+m\d)e^{2\pi\I mz}\vp\nn_m(s;q)=0.
\end{equation}
 By the mean-value theorem, there exists a number $0<\t<1$
 such that  
\begin{equation*}
\label{for:as vp}
 |\vp\nn_m(s;q)|
\le \int^{1}_{0}t^{\Re(s)-\n-1}e^{\frac{t}{1-q}}dt
={\t}^{\Re(s)-\n-1}e^{\frac{\t}{1-q}}.
\end{equation*}
 Using the Stirling formula \eqref{for:stirling} again, we obtain
 \eqref{for:lim-main}.    
\end{proof}

%===========================================================================
\section{Concluding remarks}
\label{sec:cong}
%===========================================================================

 Let $\z_q(s,t):=\z_q(s,t,1)=\sum^{\infty}_{n=1}q^{nt}{[n]_q}^{-s}$ and 
 $\z\nn_q(s):=\z_q(s,s-\n)$ for $\n\in\bN$. As a final remark, we show 
 the $q$-analogue $\z\nn_q(s)$ of $\z(s)$ can be obtained from the
 function $Z_q(s,t)$ defined by 
\begin{equation}
\label{def:Z_q}
 Z_q(s,t):=\prod^{\infty}_{l=0}(1-q^{t+l})^{-\binom{s+l-1}{l}.}
\end{equation}
 Assume $\Re(t)>0$. Then 
\begin{align}
\label{for:log-Z_q}
 \log{Z_q(s,t)}
&=-\sum^{\infty}_{l=0}\binom{s+l-1}{l}\log{(1-q^{t+l})}
=\sum^{\infty}_{l=0}\binom{s+l-1}{l}\sum^{\infty}_{n=1}\frac{q^{n(t+l)}}{n}\nonumber\\
&=\sum^{\infty}_{n=1}\frac{q^{nt}}{n}\sum^{\infty}_{l=0}\binom{s+l-1}{l}q^{nl}
=\sum^{\infty}_{n=1}\frac{q^{nt}}{n}(1-q^n)^{-s}.
\end{align}
 Since the most right hand side converges absolutely in $\Re(t)>0$, the  
 infinite product \eqref{def:Z_q} converges absolutely in this
 region. Though it is hard to expect any Euler product, the next
 proposition claims that $\z\nn_q(s)$ can be gotten by a specialization
 of the logarithmic derivative of $Z_q(s,t)$.  

%%%%%%%%%%%%%%%%%%%%%%%%%%%%%%%%%%%%%%%%%%%%%%%%%%%%%%%%%%%%%%%%%%%%%%%%%%%%
\begin{prop}
\label{prop:zZ}
 It holds that  
\begin{equation}
\label{for:zZ}
 \z\nn_q(s)
=\frac{(1-q)^s}{\log{q}}\frac{\p}{\p t}\log{Z_q(s,t)}\Bigl|_{t=s-\n}.
\end{equation}
\end{prop} 
\begin{proof}
 From \eqref{for:log-Z_q}, we have
\[
 \frac{\p}{\p t}\log{Z_q(s,t)}
=(\log{q})\sum^{\infty}_{n=1}q^{nt}(1-q^n)^{-s}
=\frac{\log{q}}{(1-q)^s}\z_q(s,t).
\]
 Hence the claim follows.
\end{proof}

%%%%%%%%%%%%%%%%%%%%%%%%%%%%%%%%%%%%%%%%%%%%%%%%%%%%%%%%%%%%%%%%%%%%%%%%%%%%
\begin{remark}
 Let $p$ be a prime number and $q=p^{-m}$ $(m\in\bN)$. Let
 $\bF_{q^{-n}}$ be the field of $q^{-n}$ elements. By
 \eqref{for:log-Z_q}, $Z_q(s,t)$ can be written also as       
\[
 Z_q(s,t)
=\exp\Bigl(\sum^{\infty}_{n=1}(q^{-n}-1)^{-s}\frac{u^n}{n}\Bigr)
=\exp\Bigl(\sum^{\infty}_{n=1}(\#\bF_{q^{-n}}^{\times})^{-s}\frac{u^n}{n}\Bigr),
\]
 where $u=q^{t-s}$. 
\end{remark}

 To see basic properties of the function $Z_q(s,t)$, we recall
 Appell's $\cO$-function $\cO_q(t;{\bsym{\om}})$ defined by   
\[
 \cO_q(t;{\bsym{\om}})
:=\prod_{l_1,\ldots,l_m\ge
0}\bigl(1-q^{l_1\om_1+\cdots+l_m\om_m+t}\bigr),
\]
 where $\bsym{\om}:=(\om_1,\ldots,\om_m)\in\bC^m$ with $\Re(\om_i)>0$
 $(1\le i\le m)$ (see \cite{Appell1882}, also
 \cite{KurokawaWakayama2004}). Similarly to the discussion in
 \cite{HashimotoWakayama}, we have the 

%%%%%%%%%%%%%%%%%%%%%%%%%%%%%%%%%%%%%%%%%%%%%%%%%%%%%%%%%%%%%%%%%%%%%%%%%%%%
\begin{prop}
 $(i)$\ We have $Z_q(0,t)=(1-q^t)^{-1}$ and
\begin{align}
\label{for:sv Z_q(m,t)}
 Z_q(m,t)=\cO_q(t;{\bf 1}_m)^{-1}, \qquad 
 Z_q(-m,t)=\prod^{m}_{l=0}(1-q^{t+l})^{(-1)^{l-1}\binom{m}{l}}
 \quad (m\in\bN),
\end{align} 
 where ${\bf 1}_m=(\underbrace{1,\ldots,1}_{m})$.

 $(ii)$\ For $m\in\bN$, we have
\begin{align}
\label{for:recursions}
\begin{array}{l}
 \displaystyle{Z_q(s+m,t)=\prod^{\infty}_{l=0}Z_q(s,t+l)^{\binom{m+l-1}{l}},\quad\,
 Z_q(s-m,t)=\prod^{m}_{l=0}Z_q(s,t+l)^{(-1)^l\binom{m}{l}}},\\[14pt]
 \displaystyle{Z_q(s,t+m)=\prod^{m}_{l=0}Z_q(s-l,t)^{(-1)^l\binom{m}{l}},\quad
 Z_q(s,t-m)=\prod^{\infty}_{l=0}Z_q(s-l,t)^{\binom{m+l-1}{l}}}.
\end{array}
\end{align} 
 In particular,
\begin{equation}
\label{for:rec}
 Z_q(s,t)=Z_q(s-1,t)Z_q(s,t+1).
\end{equation}
\end{prop}
\begin{proof}
 Equations \eqref{for:sv Z_q(m,t)} are obvious from the definition. One
 can check the ladder relations in \eqref{for:recursions} directly. 
\end{proof}

 Formulas \eqref{for:recursions} together with \Prop{zZ} yield the
 following relation for $\z\nn_q(s)$. 

%%%%%%%%%%%%%%%%%%%%%%%%%%%%%%%%%%%%%%%%%%%%%%%%%%%%%%%%%%%%%%%%%%%%%%%%%%%%%
\begin{cor}
 For $m\in\bZ$ satisfying $1\le m \le \n-1$, we have
\begin{equation}
\label{for:recm}
 \z\tp{\n-m}_q(s)=\sum^{m}_{l=0}(-1)^l\binom{m}{l}(1-q)^l\z\tp{\n-l}_q(s-l). 
\end{equation} 
 In particular,
\begin{equation}
\label{for:rec1}
 \z\tp{\n}_q(s)=\z\tp{\n-1}_q(s)+(1-q)\z\tp{\n-1}_q(s-1). 
\end{equation}
\end{cor}
\begin{proof}
 The equation \eqref{for:recm} follows from \eqref{for:zZ} and
 \eqref{for:recursions}.
\end{proof}

%%%%%%%%%%%%%%%%%%%%%%%%%%%%%%%%%%%%%%%%%%%%%%%%%%%%%%%%%%%%%%%%%%%%%%%%%%%%%
\begin{remark}
 The equation \eqref{for:rec1} was found in \cite{Zhao}. One can show it
 directly from the definition. More generally, we have
\[
 \z\tp{\n}_q(s,z)=\z\tp{\n-1}_q(s,z)+(1-q)\z\tp{\n-1}_q(s-1,z).
\]
\end{remark}

%===========================================================================
% Reference
%===========================================================================

\smallskip

%===========================================================================
% Author
%===========================================================================

\textsc{Masato Wakayama}\\
Faculty of Mathematics, Kyushu University,
Hakozaki Fukuoka 812-8581, Japan.\\
e-mail : \texttt{wakayama@math.kyushu-u.ac.jp}\\

\textsc{Yoshinori Yamasaki}\\
Graduate School of Mathematics, Kyushu University,
Hakozaki Fukuoka 812-8581, Japan.\\
e-mail : \texttt{ma203032@math.kyushu-u.ac.jp}

%===========================================================================
%===========================================================================

\begin{thebibliography}{GaWar}
%%%%%%%%%%%%%%%%%%%%%%%%%%%%%%%%%%%%%%%%%%%%%%%%%%%%%%%%%%%%%%%%%%%%%%%%%%%%%
\bibitem[1]{AndrewsAskeyRoy1999}
 G.E. Andrews, R. Askey and R. Roy (1999)
 Special Functions.
 Encyclopedia of Mathematics and its Applications, vol. 71.
 Cambridge University Press. 
%%%%%%%%%%%%%%%%%%%%%%%%%%%%%%%%%%%%%%%%%%%%%%%%%%%%%%%%%%%%%%%%%%%%%%%%%%%%%
\bibitem[2]{Appell1882}
 P. Appell (1882)
 Sur une classe de fonctions analogues aux fonctions Eul\'eriennes.
 Math Ann {\bf 19}: 84--102. 
%%%%%%%%%%%%%%%%%%%%%%%%%%%%%%%%%%%%%%%%%%%%%%%%%%%%%%%%%%%%%%%%%%%%%%%%%%%%%
\bibitem[3]{KanekoKurokawaWakayama2003}
 M. Kaneko, N. Kurokawa and M. Wakayama (2003)
 A variation of Euler's approach to values of the Riemann zeta function.
 Kyushu J Math {\bf 57}: 175--192.
%%%%%%%%%%%%%%%%%%%%%%%%%%%%%%%%%%%%%%%%%%%%%%%%%%%%%%%%%%%%%%%%%%%%%%%%%%%%%
\bibitem[4]{KawagoeWakayamaYamasaki}
 K. Kawagoe, M. Wakayama and Y. Yamasaki
 $q$-Analogues of the Riemann zeta, the Dirichlet $L$-functions, 
 and a crystal zeta function.
 preprint 2004, math.NT/0402135.
%%%%%%%%%%%%%%%%%%%%%%%%%%%%%%%%%%%%%%%%%%%%%%%%%%%%%%%%%%%%%%%%%%%%%%%%%%%%%
\bibitem[5]{Koblitz1982}
 N. Koblitz (1982)
 On Carlitz's $q$-Bernoulli numbers.
 J Number Theory {\bf 14}: 332--339. 
%%%%%%%%%%%%%%%%%%%%%%%%%%%%%%%%%%%%%%%%%%%%%%%%%%%%%%%%%%%%%%%%%%%%%%%%%%%%%
\bibitem[6]{KurokawaWakayama2004}
 N. Kurokawa and M. Wakayama (2004)
 Absolute tensor products.
 Int Math Res Not {\bf 5}: 249--260.
%%%%%%%%%%%%%%%%%%%%%%%%%%%%%%%%%%%%%%%%%%%%%%%%%%%%%%%%%%%%%%%%%%%%%%%%%%%%%%
\bibitem[7]{HashimotoWakayama}
 Y. Hashimoto and M. Wakayama
 Hierarchy of the Selberg zeta functions.
 preprint 2004, math.ph/0501047.
%%%%%%%%%%%%%%%%%%%%%%%%%%%%%%%%%%%%%%%%%%%%%%%%%%%%%%%%%%%%%%%%%%%%%%%%%%%%%
\bibitem[8]{Riemann1859}
 B. Riemann (1859)
 \"Uber die Anzahl der Primzahlen unter eine gegebener Gr\"osse.
 Monatsberichte Akad Berlin: 671--680. 
%%%%%%%%%%%%%%%%%%%%%%%%%%%%%%%%%%%%%%%%%%%%%%%%%%%%%%%%%%%%%%%%%%%%%%%%%%%%%%
%\bibitem[Ti]{Titchmarsh1986}
% E.C. Titchmarsh.:
% The Theory of the Riemann Zeta-function, 2nd edn,
% Oxford Science Publications, 1986.
%%%%%%%%%%%%%%%%%%%%%%%%%%%%%%%%%%%%%%%%%%%%%%%%%%%%%%%%%%%%%%%%%%%%%%%%%%%%%
\bibitem[9]{Satoh1989}
 J. Satoh (1989)
 $q$-analogue of Riemann's $\z$-function and $q$-Euler numbers.
 J Number Theory {\bf 31}: 346--362.
%%%%%%%%%%%%%%%%%%%%%%%%%%%%%%%%%%%%%%%%%%%%%%%%%%%%%%%%%%%%%%%%%%%%%%%%%%%%%%
\bibitem[10]{Tsumura1999}
 H. Tsumura (1999)
 A note on $q$-analogues of Dirichlet series.
 Proc Japan Acad Ser A Math Sci {\bf 75}: 23--25.
%%%%%%%%%%%%%%%%%%%%%%%%%%%%%%%%%%%%%%%%%%%%%%%%%%%%%%%%%%%%%%%%%%%%%%%%%%%%%%
\bibitem[11]{Tsumura2001}
 H. Tsumura (2001)
 On modification of the $q$-$L$-series and its applications.
 Nagoya Math J {\bf 164}: 185--197.
%%%%%%%%%%%%%%%%%%%%%%%%%%%%%%%%%%%%%%%%%%%%%%%%%%%%%%%%%%%%%%%%%%%%%%%%%%%%%%
\bibitem[12]{WhittakerWatson1927}
 E.T. Whittaker and G.N. Watson (1927)
 A Course of Modern Analysis. 4th edn.
 Cambridge University Press.
%%%%%%%%%%%%%%%%%%%%%%%%%%%%%%%%%%%%%%%%%%%%%%%%%%%%%%%%%%%%%%%%%%%%%%%%%%%%%
\bibitem[13]{Zhao}
 J. Zhao
 $q$-Multiple Zeta Functions and $q$-Multiple Polylogarithms.
 preprint 2003, math.QA/0304448.
%%%%%%%%%%%%%%%%%%%%%%%%%%%%%%%%%%%%%%%%%%%%%%%%%%%%%%%%%%%%%%%%%%%%%%%%%%%%%
%%%%%%%%%%%%%%%%%%%%%%%%%%%%%%%%%%%%%%%%%%%%%%%%%%%%%%%%%%%%%%%%%%%%%%%%%%%%%%
%\bibitem[AAR]{AndrewsAskeyRoy1999}
% G.E. Andrews, R. Askey and R. Roy.:
% Special Functions, 
% Encyclopedia of math. and appl. vol. 71, Cambridge U.P., 1999. 
%%%%%%%%%%%%%%%%%%%%%%%%%%%%%%%%%%%%%%%%%%%%%%%%%%%%%%%%%%%%%%%%%%%%%%%%%%%%%%
%\bibitem[A]{Appell1882}
% P. Appell.:
% Sur une classe de fonctions analogues aux fonctions Eul\'eriennes,
% {\it Math. Ann.} {\bf 19} (1882), 84--102. 
%%%%%%%%%%%%%%%%%%%%%%%%%%%%%%%%%%%%%%%%%%%%%%%%%%%%%%%%%%%%%%%%%%%%%%%%%%%%%%
%\bibitem[KKW]{KanekoKurokawaWakayama2003}
% M. Kaneko, N. Kurokawa and M. Wakayama.:
% A variation of Euler's approach to values of the Riemann zeta function,
% {\it Kyushu J. Math.} {\bf 57} (2003), 175--192.
%%%%%%%%%%%%%%%%%%%%%%%%%%%%%%%%%%%%%%%%%%%%%%%%%%%%%%%%%%%%%%%%%%%%%%%%%%%%%%
%\bibitem[KWY]{KawagoeWakayamaYamasaki}
% K. Kawagoe, M. Wakayama and Y. Yamasaki.:
% $q$-Analogues of the Riemann zeta, the Dirichlet $L$-functions, 
% and a crystal zeta function,
% preprint 2004, math.NT/0402135.
%%%%%%%%%%%%%%%%%%%%%%%%%%%%%%%%%%%%%%%%%%%%%%%%%%%%%%%%%%%%%%%%%%%%%%%%%%%%%%
%\bibitem[Ko]{Koblitz1982}
% N. Koblitz.:
% On Carlitz's $q$-Bernoulli numbers,
% {\it J. Number Theory} {\bf 14} (1982) 332--339. 
%%%%%%%%%%%%%%%%%%%%%%%%%%%%%%%%%%%%%%%%%%%%%%%%%%%%%%%%%%%%%%%%%%%%%%%%%%%%%%
%\bibitem[KW]{KurokawaWakayama2004}
% N. Kurokawa and M. Wakayama.:
% Absolute tensor products,
% {\it Int. Math. Res. Not.} {\bf 5} (2004), 249--260.
%%%%%%%%%%%%%%%%%%%%%%%%%%%%%%%%%%%%%%%%%%%%%%%%%%%%%%%%%%%%%%%%%%%%%%%%%%%%%%%
%\bibitem[HW]{HashimotoWakayama}
% Y. Hashimoto and M. Wakayama.:
% Hierarchy of the Selberg zeta functions,
% preprint 2005, math.ph/0501047.
%%%%%%%%%%%%%%%%%%%%%%%%%%%%%%%%%%%%%%%%%%%%%%%%%%%%%%%%%%%%%%%%%%%%%%%%%%%%%%
%\bibitem[R]{Riemann1859}
% B. Riemann.:
% \"Uber die Anzahl der Primzahlen unter eine gegebener Gr\"osse,
% {\it Monatsberichte Akad. Berlin} (1859), 671--680. 
%%%%%%%%%%%%%%%%%%%%%%%%%%%%%%%%%%%%%%%%%%%%%%%%%%%%%%%%%%%%%%%%%%%%%%%%%%%%%%%
%%\bibitem[Ti]{Titchmarsh1986}
%% E.C. Titchmarsh.:
%% The Theory of the Riemann Zeta-function, 2nd edn,
%% Oxford Science Publications, 1986.
%%%%%%%%%%%%%%%%%%%%%%%%%%%%%%%%%%%%%%%%%%%%%%%%%%%%%%%%%%%%%%%%%%%%%%%%%%%%%%
%\bibitem[S]{Satoh1989}
% J. Satoh.:
% $q$-analogue of Riemann's $\z$-function and $q$-Euler numbers,
% {\it J. Number Theory} {\bf 31} (1989), 346--362.
%%%%%%%%%%%%%%%%%%%%%%%%%%%%%%%%%%%%%%%%%%%%%%%%%%%%%%%%%%%%%%%%%%%%%%%%%%%%%%%
%\bibitem[T\,1]{Tsumura1999}
% H. Tsumura.:
% A note on $q$-analogues of Dirichlet series,
% {\it Proc. Japan Acad. Ser. A Math. Sci.} {\bf 75} (1999), 23--25.
%%%%%%%%%%%%%%%%%%%%%%%%%%%%%%%%%%%%%%%%%%%%%%%%%%%%%%%%%%%%%%%%%%%%%%%%%%%%%%%
%\bibitem[T\,2]{Tsumura2001}
% H. Tsumura.:
% On modification of the $q$-$L$-series and its applications,
% {\it Nagoya Math. J.} {\bf 164} (2001), 185--197.
%%%%%%%%%%%%%%%%%%%%%%%%%%%%%%%%%%%%%%%%%%%%%%%%%%%%%%%%%%%%%%%%%%%%%%%%%%%%%%%
%\bibitem[WW]{WhittakerWatson1927}
% E.T. Whittaker and G.N. Watson.:
% A Course of Modern Analysis, 4th edn,
% Cambridge University Press, 1927.
%%%%%%%%%%%%%%%%%%%%%%%%%%%%%%%%%%%%%%%%%%%%%%%%%%%%%%%%%%%%%%%%%%%%%%%%%%%%%%
%\bibitem[Z]{Zhao}
% J. Zhao.:
% $q$-Multiple Zeta Functions and $q$-Multiple Polylogarithms,
% preprint 2003, math.QA/0304448.
%%%%%%%%%%%%%%%%%%%%%%%%%%%%%%%%%%%%%%%%%%%%%%%%%%%%%%%%%%%%%%%%%%%%%%%%%%%%%%
\end{thebibliography}
\end{document}